\documentclass[11pt]{article}

\newcommand{\bra}{\langle}
\newcommand{\ket}{\rangle}
\newcommand{\tl}[1]{\tilde{#1}}

\newcommand{\beq}{\begin{equation}}
\newcommand{\eeq}{\end{equation}}
\newcommand{\beqs}{\begin{eqnarray}}
\newcommand{\eeqs}{\end{eqnarray}}

\def\al{\alpha}		\def\g{\gamma}		\def\del{\delta}	
\def\eps{\epsilon}  \def\la{\lambda}	\def\tht{\theta}	\def\sig{\sigma}

\newcommand{\bfa}{{\bf a}}
\newcommand{\bfb}{{\bf b}}
\newcommand{\bfc}{{\bf c}}
\newcommand{\bfu}{{\bf u}}
\newcommand{\bfv}{{\bf v}}
\newcommand{\bfx}{{\bf x}}

\newcommand{\mR}{\mathbb{R}}
\newcommand{\mZ}{\mathbb{Z}}
\newcommand{\mC}{\mathbb{C}}
\newcommand{\mH}{\mathbb{H}}
\newcommand{\mO}{\mathbb{O}}

\usepackage{amsmath,amssymb}
\usepackage{txfonts}

\usepackage[pdftex]{hyperref,color,graphicx}
\usepackage{graphicx}
\usepackage[all,cmtip]{xy}

\newcount\colveccount  
\newcommand*\colvec[1]{\global\colveccount#1  \begin{pmatrix} \colvecnext} \def\colvecnext#1{#1 \global\advance\colveccount-1
        \ifnum\colveccount>0 \\ \expandafter\colvecnext
        \else \end{pmatrix} \fi}

\begin{document}


\title{
\Large  Algebra and geometry of Hamilton's quaternions \\ \vspace{.2cm}
\large `Well, Papa, can you multiply triplets?'}
\author{{\sc Govind S. Krishnaswami and Sonakshi Sachdev}
\\ \\ \small
Chennai Mathematical Institute,  SIPCOT IT Park, Siruseri 603103, India \\ \small
 Email: {\tt govind@cmi.ac.in, sonakshi@cmi.ac.in}}

\date{\small June 10, 2016 \\
\vspace{.3cm}
Published in {\it Resonance-Journal of Science Education}, Vol. 21, Issue 6, June 2016, pp. 529-544.}
\maketitle

\begin{abstract} \normalsize

Inspired by the relation between the algebra of complex numbers and plane geometry, William Rowan Hamilton sought an algebra of triples for application to three dimensional geometry. Unable to multiply and divide triples, he invented a non-commutative division algebra of quadruples, in what he considered his most significant work, generalizing the real and complex number systems. We give a motivated introduction to quaternions and discuss how they are related to Pauli matrices, rotations in three dimensions, the three sphere, the group SU$(2)$ and the celebrated Hopf fibrations.

\end{abstract}


\normalsize

\section{Introduction}

{\it Every morning in the early part of the above-cited month [October 1843], on my coming down to breakfast, your (then) little brother William Edwin, and yourself, used to ask me, ``Well, Papa, can you multiply triplets''? Whereto I was always obliged to reply, with a sad shake of the head: ``No, I can only add and subtract them.''} W R Hamilton in a letter dated August 5, 1865 to his son A H Hamilton \cite{RPGraves}.

Imaginary and complex numbers arose in looking for `impossible' solutions to polynomial equations such as $x^2 + 1 = 0$. Even a skeptic might be convinced of their utility upon seeing the remarkable Ferro-Tartaglia-Cardano-Bombelli formula for a {\it real} root of cubic equations such as $x^3 = 20 x + 25$, expressed as sums of cube-roots of complex numbers! No such general expression in terms of real radicals is possible. Some of the mystery surrounding complex numbers was removed once they could be visualized as vectors in the plane, following the work of Euler, Wessel, Argand and Gauss. Hamilton entered the scene in 1830 by defining complex numbers $z = x + y i$ as ordered pairs (or `couples') $(x,y)$ of real numbers. He abstracted rules for adding, multiplying and dividing ordered pairs obeying the commutative, associative and distributive laws of arithmetic\footnote{Commutative means $z_1 z_2 = z_2 z_1$, associative means $z_1 (z_2 z_3) = (z_1 z_2) z_3$ while $z_1(z_2 + z_3) = z_1 z_2 + z_1 z_3$ is distributivity.}. Thus the complex numbers $\mC$ constitute a `field' just like the reals $\mR$. The manner in which the geometry of the plane was encoded in the algebra of complex numbers\footnote{The parallelogram law for adding vectors is just the sum of complex numbers, a rotation corresponds to multiplication by a complex number of unit magnitude, the cosine of the angle between vectors is given in terms of the scalar product of the corresponding complex numbers, etc.} made a deep impression on Hamilton. The problem of finding an algebra of triples $(\al, \beta, \g)$ to describe the geometry of vectors in three dimensional (3D) space haunted him for at least fifteen years. While he knew how to add and subtract triples, the problem of multiplying them and dividing by non-zero triples seemed insurmountable. What is more, by analogy with the complex numbers, he wanted the Euclidean length (square-root of the sum of squares of the components) of the product of a pair of triples to equal the product of their lengths\footnote{The cross product of vectors had not yet been discovered. Besides being non-associative $\bfa \times (\bfb \times \bfc) \ne (\bfa \times \bfb) \times \bfc$, it does {\it not} satisfy the condition $| \bfa \times \bfb| = |\bfa||\bfb|$. Cross products were invented by Grassmann within a year of Hamilton's discovery of quaternions.}. In current language, Hamilton was looking for a real, three-dimensional, normed, associative, division algebra; we now know that such an algebra does not exist. 

In 1843, Hamilton found an ingenious way around this problem, combining two major innovations. The first was to drop the commutative law for multiplication. He wrote a triple as $t = \al + \beta i + \g j$ and assumed by analogy with the complex numbers, that $i^2 = j^2 = -1$. It follows that
	\beq
	t^2 = \al^2 - \beta^2 - \g^2 + 2 \al \beta i + 2 \al \g j + \beta \g (ij + ji).
	\eeq
Moreover, he wanted the `law of moduli' $|t|^2 = |t^2|$ to hold. Now $|t|^2 = \al^2 + \beta^2 + \g^2$ would equal $|t^2|$ if $ij + ji = 0$, since
	\beq
	(\al^2 - \beta^2 - \g^2)^2 + (2 \al \beta)^2 + (2 \al \g)^2 = (\al^2 + \beta^2 + \g^2)^2.
	\eeq
Hamilton was ``tempted for a moment to fancy that'' $ij = -ji = 0$\footnote{The quotations in this paragraph are from a letter dated October 17, 1843 from Hamilton to J T Graves, which was reprinted in Vol xxv, pp.489-495 of \cite{Hamilton-on-quat}.}. But this ``seemed uncomfortable'' and he settled on the ``less harsh'' condition $ij = -ji = k$, reserving to himself to ``inquire whether $k$ was $0$ or not''. The second innovation, was to ``admit in some sense a fourth dimension of space for the purpose of calculating with triples''\footnote{In taking this step, Hamilton was partly motivated by speculations on how time and space may be ``girdled'' together, as well as some vague notions (based on Kantian philosophy) of geometry and arithmetic being the sciences of space and time respectively \cite{ETBell}.}. In other words, he assumed that $k$ was linearly  independent of $1$, $i$ and $j$. Applying the associative law, he inferred that $k^2 = ij ij = -1$, $ki = j$ and $jk = i$. Hamilton vividly describes his elation upon discovering\footnote{The article by N Mukunda in this issue of Resonance places Hamilton's discovery of quaternions in the wider context of his life and work. See also the article by R Nityananda devoted to Hamilton's work in optics.} these multiplication rules: {\it ``They started into life, or light, full grown, on the 16th of October, 1843, as I was walking with Lady Hamilton to Dublin, and came up to Brougham Bridge, which my boys have since called the Quaternion Bridge. That is to say, I then and there felt the galvanic circuit of thought closed, and the sparks which fell from it were the fundamental equations between $i,j,k$ {\it exactly such} as I have used them ever since.''} Extract from a letter dated October 15, 1858 from Hamilton to P G Tait \cite{RPGraves}.

\begin{figure}[ht] 
\center
\includegraphics[width=6cm]{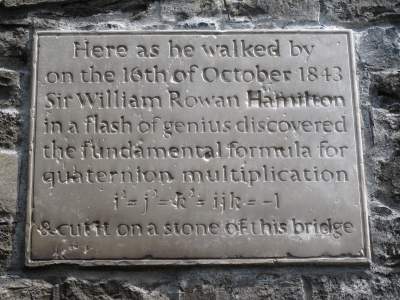}
\caption{Plaque on Broome Bridge, Dublin, commemorating Hamilton's discovery of quaternions. Photograph by Brian Dolan, Ingenious Ireland (\url{http://ingeniousireland.ie}), reprinted with permission.}
\end{figure}

In more recent language \cite{LBrand,Numbers-GTM123}, quaternions are quadruples of real numbers $q = (\al,\beta,\g,\delta)$ forming a four dimensional real vector space $\mH$ (named after Hamilton). In terms of the canonical basis
	\beq
	e = (1,0,0,0), \quad i = (0,1,0,0), \quad j = (0,0,1,0) \quad \text{and} \quad k = (0,0,0,1),
	\eeq
$q = \al e + \beta i + \g j + \delta k$. It is natural to write $q = \al e + \bfu$ where $\al e$ is called the `real' or `scalar' part of $q$. The other part $\bfu = \beta i + \g j + \delta k$ may be regarded as a vector in three-dimensional space and is called  the `imaginary' or `vectorial' part of $q$. In fact, the word {\it vector} first appears in Hamilton's writings of 1845.

Quaternions derive their power from the fact that they can be multiplied with each other, so that they form a real algebra (a vector space where vectors can be multiplied to get other vectors). Since any quaternion is a linear combination of the basis quaternions $e,i,j$ and $k$, it suffices to specify their products. First, $e$ is taken to be the multiplicative identity (sometimes denoted $1$), so $e q = qe = q$ for any quaternion $q$. In addition, $i^2 = j^2 = k^2 = -e$ while
	\beq
	ij = - ji = k, \quad 
	jk = - kj = i \quad \text{and} \quad
	ki = - ik = j.
	\eeq
Notice the analogy with cross products of unit vectors $\hat i \times \hat j = \hat k$ in 3D Euclidean space. 
\begin{table}[ht]
\begin{center}
\begin{tabular}{ |c|c|c|c|c| } 
 \hline
 $\nearrow$ & $e$ & $i$ & $j$ & $k$\\
 \hline
 $e$ & $e$ & $i$ & $j$ & $k$ \\
 \hline
 $i$ & $i$ & $-e$ & $k$ & $-j$ \\
 \hline
 $j$ & $j$ & $-k$ & $-e$ & $i$ \\
 \hline
 $k$ & $k$ & $j$ & $-i$ & $-e$ \\
 \hline
\end{tabular}
\end{center}
\caption{Hamilton's quaternion relations.}
\label{t:Hamilton-quaternion-mult}
\end{table}
The formula for the product of two quaternions $q = (\al, \bfu = \beta i + \g j + \del k)$ and $q' = (\al',\bfu' = \beta'i + \g'j + \del' k)$ can be computed using the quaternion multiplication table and written compactly in terms of the dot and cross products of their vectorial parts:
	\beqs
	q q' &=& (\al \al' + \beta\beta' i^2 + \g\g' j^2 + \del\del'k^2, ( \beta \al' + \al \beta' )i + (\g \al' + \al \g')j + (\del \al' + \al \del')k + \beta\g'(ij) \cr
	&& + \, \beta\del'(ik) + \g\beta'(ji) + \g\del'(jk) + \del\beta'(ki) + \del\g'(kj) ) \cr
	&=& \left(\al \al' - \bfu \cdot \bfu', \al \bfu' + \al' \bfu + (\beta\g' - \g\beta')k + (\del\beta' - \beta\del')j + (\g\del ' - \del\g')i \right) \cr
	&=& \left(\al \al' - \bfu \cdot \bfu', \al \bfu' + \al' \bfu + \bfu \times \bfu' \right).
	\label{e:quat-multiplication-vect}
	\eeqs
Two quaternions commute $q q' = q' q$ iff their vector parts are collinear: $\bfu \times \bfu' = 0$. Moreover, the two-sided multiplicative inverse (reciprocal) of $q$ is $q^{-1} = (\al^2 + \bfu^2)^{-1} (\al, - \bfu)$. Quaternion multiplication allows us to define polynomials in a quaternion variable. Interestingly, a polynomial of the $n^{\rm th}$ degree may have more than $n$ quaternion roots, see Box 1 for an example.

{\noindent \bf Conjugation and norm:} The conjugate of a complex number $z = x + yi$ is $\bar z = x - y i$. It allows us to define the scalar product of a pair of complex numbers $\bra w, z \ket = \Re (w \bar z) = ux + v y$ where $w = u + v i$ and $\Re$ denotes the real part. This leads to the familiar squared-length of a complex number $|z|^2 = \bra z, z \ket = \Re (z \bar z) = x^2 + y^2$. The reciprocal (left and right multiplicative inverse) of $z$ is then $\bar z/|z|^2$. By analogy with the field of complex numbers $\mC$, the conjugate of a quaternion $q = \al e + \beta i + \g j + \delta k$ is defined as $\bar q = \al e - \beta i - \g j - \del k$. One checks that the conjugate of $qq'$ is $\bar q' \bar q$. The scalar product of $q$ with another quaternion $q' = (\al', \bfu')$,
	\beq
	\bra q, q' \ket = \Re (q \bar q') = \Re (\bar q q') = \al \al' + \bfu \cdot \bfu' = \al \al' + \beta \beta' + \g \g' + \del \del',
	\eeq
is defined so as to recover the  usual Euclidean\footnote{Though Hamilton speculated on how time and space may be combined, he did not live to see their synthesis in the Minkowski space-time of special relativity. Notice however, that the bilinear form $\Re(qq') = \al \al' - \bfu \cdot \bfu'$ on $\mH$ is the Lorentzian inner product of special relativity.} squared-norm $ |q|^2 = \Re (q \bar q) = \bra q, q \ket = \al^2 + \beta^2 + \g^2 + \del^2$.

{\noindent \bf Law of moduli and Euler's four-square identity:} As noted before, Hamilton was led to his quaternion product by requiring that the norm of the product of two quaternions is the product of their norms, as is the case for complex numbers. For a pair of complex numbers, $|zw| = |z||w|$ follows from the `two-square identity':
	\beq
	(x^2 + y^2)(u^2 + v^2) = (xu - yv)^2 + (xv + yu)^2,
	\eeq
which goes back to Diophantus and Brahmagupta. The analogous statement $|qq'| = |q||q'|$ for two quaternions requires the following identity:
	\beqs
	(\al^2 + \beta^2 + \g^2 + \del^2)(\al'^2 + \beta'^2 + \g'^2 + \del'^2) &=& (\al\al' - \beta\beta' - \g \g' - \del\del')^2 + (\al \beta' + \beta \al' + \g \del' - \del \g')^2 \cr 
	&& + (\al \g' + \g \al' + \del \beta' - \beta \del')^2 + (\al \del' + \del \al' + \beta \g' - \g \beta')^2.
	\eeqs
Remarkably, Euler had discovered this identity in 1749, while trying to prove Fermat's conjecture (1659) that every natural number is a sum of four squares. The identity also appears in an unpublished note of Gauss from 1819.


\begin{center}
\fbox{\begin{minipage}{39em}
{\bf Box 1: Quaternion roots of polynomial equations in one variable:} The number of real roots of a real  polynomial of the $n^{\rm th}$ degree (say $x^2 + 1$, $x^3 + 1$ or $x^2 - 1$) can be anywhere between $0$ and $n$. The fundamental theorem of algebra guarantees that a polynomial of degree $n$ with complex coefficients has precisely $n$ complex roots (counting multiplicity). On the other hand, a polynomial of degree $n$ could have more than $n$ (even infinitely many!) quaternion roots. For instance, $i,j$ and $k$ are obvious roots of $q^2 + 1$. In fact, since $q^2 = (\al e + \bfu)^2 = \al^2 - \bfu \cdot \bfu + 2 \al \bfu$, every purely imaginary quaternion $q = (0,\bfu)$ of unit magnitude $|q|^2 = \bfu \cdot \bfu = 1$ is a root, and solutions may be identified with points on the two dimensional sphere $S^2$.
\end{minipage}}
\end{center}

\section{Quaternions represented as $2 \times 2$ matrices}

The basis quaternions $i,j$ and $k$ anti-commute ($ij = -ji$ etc.) and square to $-e$. For those who have met them, this is reminiscent of the Pauli matrices  arising in the quantum mechanics of spin:
	\beq
	\sig_1 = \colvec{2}{0 & 1}{1 & 0}, \quad
	\sig_2 = \colvec{2}{0 & -i}{i & 0} \quad \text{and} \quad
	\sig_3 = \colvec{2}{1 & 0}{0 & -1}.
	\eeq
Indeed, the Pauli matrices anti-commute ($\sig_1 \sig_2 = - \sig_2 \sig_1 = i \sig_3$ and cyclic permutations) and square to the identity matrix $I$. If we can get the signs right, it seems plausible that the quaternion algebra can be represented in terms of $2 \times 2$ {\it complex} matrices.

Let us first see how this is done for complex numbers, which may be represented by certain anti-symmetric $2 \times 2$ {\it real} matrices
	\beq
	z = x + y i \mapsto A(z) = \colvec{2}{x & - y}{y & x} = x I - y (i\sig_2).
	\label{e:map-complex-to-anti-symm-mat}
	\eeq
The set of matrices of this sort is an algebra: it is closed under matrix addition and multiplication. In fact, the map $A$ is an isomorphism of real algebras: $A(zz') = A(z) A(z')$, $A(z+ z') = A(z) + A(z')$ and $A(\la z) = \la A(z)$ for $\la \in \mR$. Remarkably, though multiplication of matrices is generally non-commutative, matrices of the above sort commute with each other (as (\ref{e:map-complex-to-anti-symm-mat}) involves only one of the Pauli matrices $\sig_2$), simulating the commutative multiplication of complex numbers. $z = 1 + 0 i$ is represented by the identity matrix and $|z|^2 = x^2 + y^2$ is $\det A(z)$. The reciprocal of a non-zero complex number is mapped to the inverse matrix. The map $A$ takes the set of unit complex numbers $z = e^{i \tht}$ comprising the $1$-sphere $S^1: x^2 + y^2 = 1$, to the space SO$(2)$ of $2 \times 2$ orthogonal\footnote{A matrix is orthogonal if its rows (or columns) are orthonormal: $R^t R = R R^t = I$, where $t$ denotes transposition. The set of $n \times n$ orthogonal matrices of unit determinant SO$(n)$ implement rotations in $n$-dimensional Euclidean space.} matrices $R_\tht = \colvec{2}{\cos \tht & - \sin \tht}{\sin \tht & \cos \tht}$ with unit determinant.

Similarly, quaternions may be represented as $2 \times 2$ {\it complex} matrices by exploiting the Pauli matrix algebra. Hamilton's relations (Table \ref{t:Hamilton-quaternion-mult}) are satisfied if we represent the basis quaternions $e,i,j$ and $k$ by the following matrices (this choice is not unique!)
	\beq
	F(e) = \colvec{2}{1 & 0}{0 & 1} = I, \quad
	F(i) = \colvec{2}{i & 0}{0 & -i} = i \sigma_3, \quad
	F(j) = \colvec{2}{0 & 1}{-1 & 0} = i \sigma_2, \quad
	F(k) = \colvec{2}{0 & i}{i & 0} = i \sigma_1.
	\eeq
This map extends by linearity to a $1$-$1$ homomorphism of $\mH$ into the associative algebra of $2 \times 2$ complex matrices:
	\beq
	q = \al e + \beta i + \g j + \del k \mapsto F(q) = \colvec{2}{\al + \beta i & \g + \del i}{- \g + \del i & \al - \beta i} = \colvec{2}{z & \bar w}{-w & \bar z}.
	\label{e:map-quaternions-to-complex-matrices}
	\eeq
In this representation, the non-commutativity of quaternion multiplication is reflected in that of the matrices. As before, the reciprocal of a non-zero quaternion is mapped to the inverse matrix and $|q|^2 = \al^2 + \beta^2 + \g^2 + \del^2$ is det$F(q)$. Under this map, the image of the set of unit quaternions ($|q| = 1$) is the space SU$(2)$ of $2 \times 2$ unitary matrices ($F F^\dag = F^\dag F = I$, where $\dag$ denotes complex conjugate transpose) with unit determinant. Notice moreover that unit quaternions comprise the 3-sphere $S^3: \al^2 + \beta^2 + \g^2 + \del^2 = 1$ embedded in $\mR^4$.

\section{Quaternions and rotations} 

Rotations are an important symmetry of many problems in geometry, physics and astronomy and it is useful to have explicit formulae for how rotations transform vectors. Hamilton developed quaternions in part to describe rotations in 3D, inspired by how complex numbers describe rotations in the plane.  Multiplication of $w = u + v i$ by a complex number $z = e^{i \tht}$ of unit magnitude represents a counter-clockwise rotation $R_\tht$ of ${\bf w} = (u,v)$ by angle $\tht$. These rotations may be composed $R_{\tht} R_{\phi} = R_{\tht + \phi}$ and undone $R_{\tht}^{-1} = R_{-\tht}$. Thus, unit complex numbers $e^{i \tht}$ form the circle {\it group} U$(1)$ of $1 \times 1$ unitary matrices. The map $A$ of Eq. (\ref{e:map-complex-to-anti-symm-mat}) restricted to $|z| = 1$ is an isomorphism from U$(1)$ to the group SO$(2)$ of rotations of vectors in a plane.

\begin{center}
\fbox{\begin{minipage}{39em}
The ability to multiply real numbers gives the unit reals $S^0 = \{-1,1\}$ the structure of the cyclic group of order two: $\mZ/2\mZ$. The multiplication of complex numbers gives the unit circle $S^1$ in the complex plane the structure of the group U$(1)$. Similarly, the existence of the quaternions endows the unit quaternions (i.e. the $3$-sphere $S^3$) with the structure of the group SU$(2)$.
\end{minipage}}
\end{center}

Quaternions give a convenient way of expressing rotations in three dimensions. First, we identify vectors in 3D space with imaginary quaternions $\Im \mH \cong \mR^3$. By analogy with the unit complex numbers acting as rotations of the plane, the unit quaternions ($|q| = 1$) act as rotations of $\Im \mH$. However, the action is not by multiplication, but by conjugation $v \mapsto v' = qvq^{-1}$, which takes $\Im \mH$ to $\Im \mH$. Indeed, one checks using (\ref{e:quat-multiplication-vect}) that if $v$ is purely imaginary, then so is $v'$. Moreover by the law of moduli $|v'| = |qvq^{-1}| =  |v|$ so that lengths of vectors are preserved under conjugation. Explicitly, if $v = (0, \bfv)$, $q = (\al, \bfu)$ and $q^{-1} = \bar q = (\al, -\bfu)$ with $\al^2 + |\bfu|^2 = 1$, then using (\ref{e:quat-multiplication-vect}) we get
	\beqs
	qvq^{-1} &=& \left( - \bfu \cdot \bfv, \al \bfv + \bfu \times \bfv \right) (\al, - \bfu) \cr
	&=& \left(0, \al^2 \bfv + 2 \al (\bfu \times \bfv) + (\bfu \cdot \bfv) \bfu + \bfu \times (\bfu \times \bfv) \right) \cr
	&=& \left(0, (\al^2 - |\bfu|^2) \bfv + 2\al (\bfu \times \bfv) + 2(\bfu \cdot \bfv)\bfu \right) \equiv (0, R \, \bfv).
	\label{e:conj-imaginary-q-by-unit-q}
	\eeqs
Here $R = R(q)$ is the $3 \times 3$ matrix $R = (\al^2 - |\bfu|^2) I - 2\al \tl u + 2 \bfu \bfu^t$ with $\tl u$ the matrix whose entries are $\tl u_{ab} = \sum_c \eps_{abc} \bfu_c$\footnote{Here $\eps_{abc}$ is the Levi-Civita symbol, anti-symmetric under exchange of any pair of indices and with $\eps_{123} = 1$.}. Thus, conjugation is an origin and length-preserving linear transformation of $\Im \mH$. Such transformations $R$ must either be proper rotations ($\det R = 1$, comprising the group SO$(3)$) or rotations composed with a reflection ($\det R = -1$). It is possible to show that the space SU$(2)$ of unit quaternions is path connected, so that the continuous function $\det R(q)$ must always equal the value $1$ that it takes at the identity quaternion $q=e$. Thus, conjugation by unit quaternions act as proper rotations. Now, Euler showed in 1776 that any rotation of 3D space may be regarded as a counter-clockwise rotation by angle $\tht$ about an axis specified by a unit vector $\hat n$. It remains to relate the parameters $\al$ and $\bfu$ of the rotation $R(q)$ to $\tht$ and $\hat n$. It turns out that $\al = \cos(\tht/2)$ and $\hat n$ is the unit vector along $\bfu$.

\begin{figure}[ht] 
\center
\includegraphics[width=4cm]{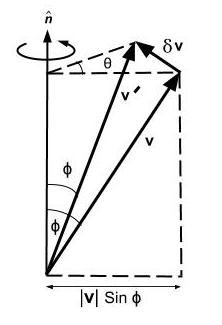}
\caption{The infinitesimal change $\del \bfv \approx \theta \: \hat n \times \bfv$ in a vector $\bfv$ due to a counter-clockwise rotation by a small angle $\theta$ about the axis $\hat n$. Notice that $|\del \bfv| \approx \tht \, |\bfv| \sin \phi$ and that $\del \bfv$ points along $\hat n \times \bfv$.}
\label{f:infinitesimal-rotn}
\end{figure}

To see this, we first consider an infinitesimal rotation, which should correspond to conjugation by a unit quaternion $q = (\al, \bfu)$ that differs infinitesimally from the identity $q = e = (1, \bf0)$ (i.e., $\al^2 \approx 1$ and $|\bfu|^2 \approx 0$). In this case, from Eq. (\ref{e:conj-imaginary-q-by-unit-q}) we find 
	\beq
	R \bfv = \bfv' = \bfv + \del \bfv  \approx \bfv + 2 \al \bfu \times \bfv.
	\label{e:infinitesimal-rotn}
	\eeq
A delightful piece of trigonometry shows that this represents a counter-clockwise rotation of $\bfv$ about the axis $\hat n = \bfu/|\bfu|$ by small angle $\tht$ with $|\tht| = 2 |\al \bfu| = 2 |\bfu| \sqrt{1 - |\bfu|^2} \approx 2 |\bfu|$ since $\al^2 = 1 - |\bfu|^2$ and $|\bfu|^2 \approx 0$ (see Fig. \ref{f:infinitesimal-rotn}). More generally, we may always express the {\it unit} quaternion $q = \al e + \bfu$ as\footnote{Recall Euler's formula $e^{i \tht} = \cos \tht + i \sin \tht$ for a unit complex number as an exponential of an imaginary number. Similarly, unit quaternions are expressible as exponentials of purely imaginary quaternions $e^{(\tht/2) {\hat n}} = \cos (\tht/2) e + \sin(\tht/2) {\hat n}$. To show this, notice that the unit imaginary quaternion ${\hat n} = n_x i + n_y j + n_z k$ satisfies ${\hat n}^2 = -e$, ${\hat n}^3 = - {\hat n}$, ${\hat n}^4 = e$ etc., allowing one to sum the exponential series.}
	\beq
	q = \cos(\tht/2) \: e + \sin(\tht/2) \: \hat n,
	\label{e:unit-q-in-terms-tht-n-hat}
	\eeq
where $\hat n$ is the unit vector along $\bfu$. The division by two ensures that for small $|\bfu|$, $\tht$ may be interpreted as the angle of rotation, since $|2 \bfu| = |2 \sin(\tht/2) \hat n| \approx |\tht|$. Plugging $\al = \cos(\tht/2)$ and $\bfu = \sin(\tht/2) \: \hat n$ into $R$ in Eq. (\ref{e:conj-imaginary-q-by-unit-q}) we arrive at
	\beqs
	R \bfv &=&  \left(\cos^2\left(\frac{\theta}{2} \right) - \sin^2\left(\frac{\theta}{2} \right)\right)\bfv  +  2 \cos\left(\frac{\theta}{2}\right)\sin\left(\frac{\theta}{2} \right) (\hat n \times \bfv) + 2\sin^2\left(\frac{\theta}{2} \right)( \hat n \cdot \bfv) \:\hat n \cr \cr
	&=& \text{cos} \theta \: \bfv +  \text{sin} \theta \: \hat n \times \bfv + (1  - \text{cos} \theta) (\hat n \cdot \bfv) \:\hat n .
	\label{e:Rodrigues-formula}
	\eeqs
Interestingly, this formula for the result of rotating a  vector $\bfv$ by angle $\theta$ about an axis $\hat n$ was obtained by B O Rodrigues (1840) shortly before Hamilton's discovery of quaternions, using Euler's 4-square identity. For small $\tht$, Eq. (\ref{e:Rodrigues-formula}) reduces to Eq. (\ref{e:infinitesimal-rotn}). In general, the rotated vector $\bfv'$ is expressed as a linear combination of (the generally non-orthogonal vectors) $\bfv$, $\hat n \times \bfv$ and $\hat n$. For planar rotations ($\bfv$ perpendicular to the axis $\hat n$), these vectors become orthogonal and the third term is absent.

The set of unit quaternions $\{ |q|=1 \}$ is closed under multiplication and inverses and forms a {\it group} isomorphic to SU$(2)$ via the map $F$ of Eq. (\ref{e:map-quaternions-to-complex-matrices}). The map $R(q)$ defines a homomorphism from unit quaternions (or SU$(2)$) onto the rotation group SO$(3)$. The homomorphism is $2:1$ since conjugation by both $q = e$ and $q = -e$ correspond to the identity rotation (the `kernel' is $\{ \pm e \}$). The appearance of the half-angle $\tht/2$ in Eq. (\ref{e:unit-q-in-terms-tht-n-hat}) is symptomatic of this. We say that SU$(2)$ is a double cover of SO$(3)$. This fact has many interesting consequences, such as the possibility for particles such as the electron to have spin one-half \cite{gottfried-yan,ecg-nm}.

The abelian group of unit complex numbers U$(1)$ has discrete subgroups consisting of the $n^{\rm th}$ roots of unity (the cyclic group $\mZ/n\mZ$ for each $n = 1,2,3, \ldots$). Similarly, the unit quaternions have discrete subgroups, see Box 2 for an interesting non-abelian example. Remarkably, quaternions also help in understanding rotations in {\it four} dimensions as discussed in Box 3.

\begin{center}
\fbox{\begin{minipage}{39em}
{\bf Box 2: Quaternion group:} The unit quaternions have an interesting discrete subgroup: the basis quaternions and their negatives $\pm e, \pm i, \pm j, \pm k$ form the quaternion group of order eight. Each element other than $\pm e$ generates a cyclic subgroup of order $4$, e.g. $\bra i \ket = \{ i, i^2 = -1, i^3 = -i, i^4 = e \}$. The quaternion group acts by left (or right) multiplication on $\mH$. For instance, $i(\al e + \beta i + \g j + \del k) = -\beta e + \al i - \del j + \g k$ may be interpreted as a simultaneous rotation by $90^\circ$ counter-clockwise in the $1-i$ and $j-k$ planes.
\end{minipage}}
\end{center}

\begin{center}
\fbox{\begin{minipage}{39em}
{\bf Box 3: Quaternions and rotations in four dimensions:}
Remarkably, quaternions also help us to understand rotations in $\mR^4$. As before, one may show that the group SU$(2)$ of unit quaternions act via left and right multiplication $v \mapsto q v q'$ as length and orientation preserving linear transformations of the 4D space of quaternions $v$. This gives a homomorphism from SU$(2) \times$SU$(2)$ onto the rotation group SO$(4)$ with kernel $\pm (I,I)$. Thus the former is a double cover of the latter. The group SO$(4)$ arises as a symmetry of the Kepler problem as well as the hydrogen atom. The two copies of SU$(2)$ are related to the conservation of angular momentum and Laplace-Runge-Lenz vectors.
\end{minipage}}
\end{center}

\section{Hopf maps for reals, complexes and quaternions}

\begin{figure}[ht] 
\center
\includegraphics[width=7cm]{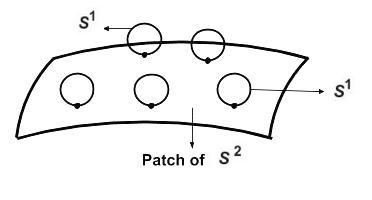}
\caption{A circle bundle over a patch of the two-sphere.}
\label{f:circle-bundle}
\end{figure}

We have met the unit circle $S^1$, the surface of the unit sphere $S^2$ and the unit reals $\{ \pm 1 \}$ denoted $S^0$. More generally, $S^n = \{ (x_1, x_2, \ldots, x_{n+1}) \in \mR^{n+1} \, | \, x_1^2 + x_2^2 + \cdots + x_{n+1}^2 = 1 \}$ is called the unit $n$-sphere. The Hopf map (H Hopf, 1931) and its generalizations are many-to-one maps between certain spheres. There are Hopf maps associated to the reals, complexes, quaternions and octonions\footnote{The octonions $\mO$ are an eight-dimensional normed division algebra generalizing $\mR, \mC$ and $\mH$. An octonion is a linear combination of the basis octonions $e_0, \cdots, e_7$, which obey a non-commutative and non-associative product rule \cite{conway-smith}.} (but no more!):
	\beq
	S^1 \to S^1, \quad S^3 \to S^2, \quad
	S^7 \to S^4 \quad \text{and} \quad S^{15} \to S^8.
	\eeq
The maps go from a `total-space' sphere `upstairs' to a `base-space' sphere `downstairs'. In each case, the pre-image of {\it any} point on the base sphere turns out (see below) to be a sphere of dimension one less than the base! These pre-images of points on the base are called fibres and the total space is called a fibre bundle (see Fig. \ref{f:circle-bundle}). The fibres in the above four cases are $S^0, S^1, S^3$ and $S^7$. This allows the sphere upstairs to be viewed locally (though not globally) as the Cartesian product of the base sphere and the fibre sphere. Thus we say that $S^3$ is a circle ($S^1$) bundle over $S^2$ and that $S^7$ is an $S^3$ bundle over $S^4$. In the first three cases, the fibres have a group structure: $S^0 \cong \mZ/2\mZ$, $S^1 \cong U(1)$ and $S^3 \cong$ SU(2), so that the base spheres $S^1$, $S^2$ and $S^4$ may be viewed as quotients of the total spaces by group actions:
	\beq
	S^1 \cong S^1/S^0, \quad 
	S^2 \cong S^3/S^1 \quad  \text{and} \quad
	S^4 \cong S^7/S^3.
	\eeq

\begin{figure}[ht] 
\center
\includegraphics[width=5cm]{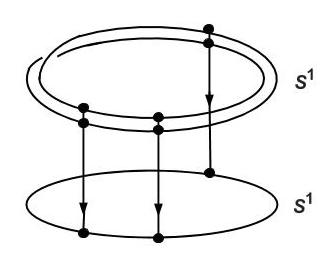}
\caption{A `cartoon' of the 2:1 real Hopf map $S^1 \to S^1$ with $S^0$ fibres. For easy visualization, the points identified `upstairs' are {\it not} antipodal.}
\label{f:Hopf-map-reals}
\end{figure}

To understand the Hopf maps in a little more detail, let us start with the real case, where it goes from the real $1$-sphere $S^1$ (which, with a view towards generalization, we denote by $S^1_{\mR}$) to the `real projective line' $P^1(\mR)$. The latter is the set of lines through the origin in $\mR^2$. The line through $(x_1,x_2) \ne (0,0)$ is an equivalence class $[(x_1,x_2)]$ of ordered pairs that differ by rescaling by a non-zero real number $\la$: $(x_1, x_2) \sim \la (x_1, x_2)$. Any such line (other than the vertical line through the origin), contains exactly one point with $x_1 = 1$, that is, a point of the form $(1,x_2)$. This point can be used to label that line. Using the label $(0,\infty)$ (or $(0,-\infty)$) for the vertical line through the origin, we see that $P^1(\mR)$ may be identified with the set of real numbers $x_2$ `compactified' by the `point' at $\pm \infty$  to form a circle, so $P^1(\mR) \cong S^1$. The Hopf map takes $(x_1, x_2) \in S^1_\mR$ to the equivalence class $[(x_1,x_2)]$. Under this map all points of the form $\la(x_1, x_2)$ with $|\la| = 1$ (i.e. where $\la = \pm 1$ is a `unit' real) are mapped to the same equivalence class $[(x_1,x_2)]$. Thus we may identify the pre-image $\{ (x_1,x_2)/|\bfx|, -(x_1, x_2)/|\bfx| \}$ of $[(x_1,x_2)]$ with $S^0 = \{1,-1 \}$ (here $|\bfx|^2 = x_1^2 + x_2^2$). The set of unit reals $S^0$ is the fibre over $[(x_1,x_2)]$. We say that the real projective line is the quotient of the real 1-sphere by the group of unit reals. This Hopf map may also be viewed as the $2:1$ map $e^{i \tht} \mapsto e^{2 i \tht}$ of the unit circle to itself, obtained by identifying antipodal points (see Fig. \ref{f:Hopf-map-reals}).

Replacing the reals with the complexes or quaternions leads to corresponding Hopf fibrations. For instance, we have the Hopf map from the complex 1-sphere
	\beq
	S^1_{\mC} = \left\{ (z_1, z_2) \in \mC^2 \: | \: |z_1|^2 + |z_2|^2 = 1 \right\} \cong S^3
	\eeq
to the complex projective line $P^1(\mC)$ taking $(z_1, z_2) \in S^1_\mC$ to the equivalence class $[(z_1, z_2)]$. As before, points of the form $\la (z_1, z_2)$ with $\la = e^{i\theta}$ a complex number of unit magnitude, are mapped to the same point/equivalence class $[(z_1, z_2)]$ in $P^1(\mC)$. So each fibre is a circle $S^1$. Furthermore, as for the reals, we may choose the equivalence class representatives as $(1,z_2)$ for $z_2 \in \mC$ and $(0,\infty)$ and thereby identify $P^1(\mC)$ with the complex plane (or $\mR^2$) compactified by the `point at infinity'. So, $P^1(\mC) \cong S^2$. Thus the Hopf map $S^3 \to S^2$ with $S^1$ fibres allows us to view $S^3$ as a circle bundle over the two-sphere (see Fig. \ref{f:circle-bundle}). Moreover, exploiting the fact that $S^3$ and $S^1$ can be identified with the groups SU$(2)$ and U$(1)$ of unit quaternions and complexes, it can be shown that $S^2$ is the quotient of SU$(2)$ by U$(1)$. This Hopf map appears in the study of polarized light \cite{rajaram}, the Dirac magnetic monopole \cite{Nakahara}, single qubit Hilbert spaces \cite{mosseri-dandoloff} as well as the planar three-body problem. Depending on the context, the `base' $S^2$ is referred to as the Bloch/Poincar\'e/shape sphere.

Turning to quaternions, we have a map from the quaternion $1$-sphere
	\beq
	S^1_{\mathbb{H}} = \{ (q_1, q_2) \in \mathbb{H}^2 \: | \: |q_1|^2 + |q_2|^2 = 1 \} \cong S^{7}
	\eeq
to the quaternion projective line $P^1(\mH) \cong S^4$ (one point compactification of $\mH \cong \mR^4$) with fibres given by the unit quaternions $S^3$. This Hopf fibration makes an appearance in the theory of instantons \cite{Nakahara} (Euclidean field configurations localized in both space and time) relevant to tunnelling phenomena in particle physics and in the study of the geometry of two qubits. The final Hopf fibration $S^7 \to S^{15} \to S^8$ results from a similar construction involving the octonions\footnote{One difference is that unlike unit reals, complexes or quaternions, the unit octonions ($S^7$) do not form a group, since the octonion multiplication is non-associative.}.

\begin{center}
\fbox{\begin{minipage}{31em}
\centerline{\bf The Hopf maps for $\mR, \mC, \mH$ and $\mO$}
\begin{center}
\leavevmode
\xymatrix{ \{ |x|=1 \} = S^0\ar[r] & S^1_\mR \cong S^1 \ar[d]\\
            & P^1(\mR)\cong S^1}
\leavevmode
\xymatrix{ \{ |z| = 1 \} \cong S^1 \ar[r] & S^1_\mC \cong S^3 \ar[d]\\
            & P^1(\mC)\cong S^2}
\end{center}
\begin{center}
\leavevmode
\xymatrix{ \{ |q| = 1 \} \cong S^3 \ar[r] & S^1_\mH \cong S^7 \ar[d]\\
            & P^1(\mH)\cong S^4}
\leavevmode
\xymatrix{  \{ |o| = 1 \} \cong S^7 \ar[r] & S^1_\mO \cong S^{15} \ar[d]\\
            & P^1(\mO)\cong S^8}
\end{center}
\end{minipage}}
\end{center}

\section{Afterword}

Hamilton thought that the discovery of quaternions in the mid $19^{\rm th}$ century was as important as the discovery of calculus at the end of the$17^{\rm th}$ century. Quaternions were vigorously researched and taught in parts of Ireland, England and the United States (a Quaternion Society, the ``International Association for Promoting the Study of Quaternions and Allied Systems of Mathematics'' was founded at Yale University in 1895). Though the importance of quaternions was probably overestimated by Hamilton and the `quaternionists', their discovery catalyzed several related developments. Hamilton's way of thinking unlocked the freedom available in relaxing one or other axiom of arithmetic, and in working with higher dimensional algebras, even to understand problems in ordinary space. Within a few months of Hamilton's discovery of quaternions, his friend J T Graves found another `hypercomplex' number system by dropping the condition of associativity: the 8-dimensional normed division algebra of octonions (rediscovered by Cayley in 1845). The contemporaneous development of the {\it non-commutative} algebra of matrices (especially beginning with the work of Cayley, recall the Cayley-Hamilton theorem) was arguably as important to the mathematical sciences and engineering as the discovery of calculus. Quaternions also fuelled the development of vector algebra and vector calculus (by Grassmann, Gibbs, Heaviside, and Helmholtz) and the formulation of the equations of particle and continuum mechanics and electromagnetism in terms of them.

Though quaternions may be viewed simply as a sub-algebra of $2 \times 2$ complex matrices, they are very special. Indeed, Frobenius (1877) showed that any finite dimensional associative division algebra over the reals, is isomorphic to one of $\mR, \mC$ or $\mH$. Relaxing the condition of associativity allows for the octonions (Hurwitz, 1923). Quaternions continue to appear in various areas of mathematical research: E.g., quaternionic versions of integers in number theory, quaternionic representations of groups and quaternionic analogues of real and complex manifolds, with many connections to physics. Quaternions are also used by computer scientists and engineers! The parametrization of rotations by unit quaternions is used by programmers to animate camera movements \cite{Shoemake} in video games (such as Tomb Raider) and to smoothly interpolate between successive orientations of an airplane in flight simulators.
 
{\noindent \bf Acknowledgements:} We thank B V Rao, R Nityananda and N Mukunda for carefully reading this article and suggesting improvements.


\end{document}